\documentclass[11pt,a4paper]{article}

\usepackage{epsf,epsfig,amsfonts,amsgen,amsmath,amstext,amsbsy,amsopn,amsthm,amsfonts,amssymb,amscd,bbding
}

\usepackage{latexsym,bm}
\usepackage{ebezier,eepic}
\usepackage{color}
\usepackage{tikz,multirow}
\usepackage{enumerate}
\usepackage{geometry}
\usepackage{indentfirst}
\usetikzlibrary{decorations.markings}
\tikzstyle{vertex}=[circle, draw, inner sep=0pt, minimum size=4pt]

\newtheorem{thm}{Theorem}[section]
\newtheorem{claim}{Claim}
\newtheorem{lemma}{Lemma}[section]

\begin{document}

\baselineskip 16pt
 %\setcounter{chapter}{1}
 % defining short form------
 \newcommand{\la}{\lambda}
 \newcommand{\si}{\sigma}
 \newcommand{\ol}{1-\lambda}
 \newcommand{\be}{\begin{equation}}
 \newcommand{\ee}{\end{equation}}
 \newcommand{\bea}{\begin{eqnarray}}
 \newcommand{\eea}{\end{eqnarray}}
\newcommand{\bL}{\b{\textit{L}}}
\newcommand{\bN}{\b{\textit{N}}}
\newcommand{\bB}{\b{\textit{B}}}

\title{\bf\Large Perfect matching and distance spectral radius in graphs and bipartite graphs}

\date{}

\author{
Yuke Zhang  ~and Huiqiu Lin\footnote{Supported
by the National Natural Science Foundation of China (Nos. 11771141 and 12011530064).\quad\quad\quad E-mail: huiqiulin@126.com (H.Q. Lin),~zhang\_yk1029@163.com (Y.K. Zhang)}\\
{\footnotesize Department of Mathematics,
East China University of Science and Technology,}\\
{\footnotesize Shanghai 200237, P.R.~China}
}
\maketitle
\vspace{-9mm}

\maketitle

\begin{abstract}
A perfect matching in a graph $G$ is a set of nonadjacent edges covering every vertex of $G$.
Motivated by recent progress on the relations between the eigenvalues and the matching number of a graph,
in this paper, we aim to present a distance spectral radius condition to guarantee the existence of a perfect matching.
Let $G$ be an $n$-vertex connected graph where $n$ is even and $\lambda_{1}(D(G))$ be the distance spectral radius of $G$.
Then the following statements are true.

\noindent$\rm{I)}$ If $4\le n\le10$ and ${\lambda }_{1} (D\left(G\right))\le {\lambda }_{1} (D(S_{n,{\frac{n}{2}}-1}))$,
then $G$ contains a perfect matching unless $G\cong S_{n,{\frac{n}{2}-1}}$ where $S_{n,{\frac{n}{2}-1}}\cong K_{{\frac{n}{2}-1}}\vee ({\frac{n}{2}+1})K_1$.

\noindent$\rm{II)}$ If $n\ge 12$ and ${\lambda }_{1} (D\left(G\right))\le {\lambda }_{1} (D(G^*))$, then $G$ contains a perfect matching unless $G\cong G^*$ where $G^*\cong K_1\vee (K_{n-3}\cup2K_1)$.

Moreover, if $G$ is a connected $2n$-vertex balanced bipartite graph with $\lambda_{1}(D(G))\le \lambda_{1}(D(B_{n-1,n-2})) $, then $G$ contains a perfect matching, unless $G\cong B_{n-1,n-2}$
where $B_{n-1,n-2}$ is obtained from $K_{n,n-2}$ by attaching two pendent vertices to a vertex in the $n$-vertex part.

 \vspace{5mm}

\vspace{3mm}

 \end{abstract}

 \baselineskip=0.30in

\section{Introduction}
All graphs considered in this paper are undirected, connected and simple.
Let $G$ be a graph with vertex set $V(G)=\{v_1,v_2,\ldots,v_n\}$
and edge set $E(G)$. The \emph{distance}  between $v_i$ and $v_j$, denoted by $d_G(v_i,\,v_j)$ (or $d_{ij}$), is the length of a shortest path from $v_i$ to $v_j$. The \emph{distance matrix} of $G$,
denoted by $D(G)$, is the $n\times n$ real symmetric matrix whose $(i,\,j)$-entry is  $d_G(v_i,\,v_j)$ $(\mbox{or }d_{ij})$, then we can order the eigenvalues of $D(G)$ as $$\lambda_1(D(G))\geq \lambda_2(D(G))\geq \cdots\geq \lambda_n(D(G)).$$  By the Perron-Frobenius theorem, $\lambda_1(D(G))$ is always positive (unless $G$ is trivial) and $\lambda_1(D(G)) \ge |\lambda_{i}(D(G))|$ for $i=2,3,\ldots, n$, and we call $\lambda_1(D(G))$ the \emph{distance spectral radius}.

The study of distance eigenvalues can be traced back to 1971 by Graham and Pollack \cite{Graham} and they described a relationship between the number of negative distance eigenvalues and the addressing problem in data communication system. In the same paper, they showed a very interesting and insightful result that the determinant of the distance matrix of a tree with order $n$ is $(-1)^{n-1}(n-1)2^{n-2}$, which is independent of the structure of the tree. Since then, the study of distance eigenvalues of a graph has become a research subject of enormous interest and this topic has received growing attention over recent years, some latest results see \cite{AH1,Lin,lhq,lzz,xj}. For more results on the distance matrix and its spectral properties, we refer the reader to the excellent survey \cite{AH}.

A \emph{matching} $M$ of a graph $G$ is a set of pairwise nonadjacent edges. The maximum number of edges of a matching in $G$ is called \emph{matching number} of $G$, denoted by $\alpha(G)$. A vertex incident with an edge in $M$ is called to be \emph{saturated} by $M$. A \emph{perfect matching} is one matching which all vertices of $G$ are saturated by it. Obviously, a graph with a perfect matching has an even number of vertices and $\alpha(G)=\frac{|V(G)|}{2}$.

The studies of the connections between eigenvalues and the matching number of a graph were mainly based on the tree in 1990s. Chang \cite{an} obtained an upper bound and a tight lower bound for the second largest eigenvalue of an $n$-vertex tree with a given matching number. Later, Hou and Li \cite{hou} gave some upper bounds for the spectral radius of a tree in terms of its order and matching number. Brouwer and Heamers \cite{Brouwer} investigated this problem to the general graph and showed that if $\mu_1(G) \le 2\mu_{n-1}(G)$, then $G$ contains a perfect matching where $\mu_1(G)$ and $\mu_{n-1}(G)$ are the largest and second smallest Laplacian eigenvalues of $G$, respectively. Feng, Yu and Zhang \cite{feng} investigated the maximal spectral radius of graphs with a given matching number and order. In the past decade, some quality and interesting results on the matching number of a graph and its distance spectral radius have been obtained. Ili\'{c} \cite{llic} characterized $ n $-vertex trees with a given matching number
which minimize the distance spectral radius. Liu \cite{liu} characterized graphs with minimum distance spectral
radius in connected graphs on $n$ vertices with fixed matching number. Zhang \cite{zhang} and Lu and Luo \cite{luluo} characterized unicyclic graphs with a perfect matching and a given
matching number which minimize the distance spectral radius, respectively.

Very recently, O \cite{suil} proved a lower bound for the spectral radius in an $ n $-vertex graph to guarantee the existence of a perfect matching. Along this line, we consider this problem with respect to the distance spectral radius in this paper.

We denote by $ G\cup H $ the \emph{disjoint union} of two graphs $G$ and $ H $, which is the graph with  $V(G\cup H)=V(G)\cup V(H)$ and $E(G\cup H)=E(G)\cup E(H)$, particularly if $G\cong H$, then write $2G= G\cup H $ for short. Denote by $G\vee H $ the \emph{join} of two graphs
$ G $ and $ H $, which is the graph such that $ V(G \vee H) = V(G) \cup V(H) $ and $ E(G \vee H) = E(G)\cup E(G) \cup \{vu:
u \in V(G)\ \mbox{and}\  v \in V(H)\}. $
The graph $S_{n,k}$ is obtained from a copy of $K_k$ by adding $n-k$ vertices, each of which has neighborhood $V(K_k)$ i.e. $S_{n,k}\cong K_k\vee(n-k)K_1$. In the following, we first give a distance spectral radius condition which guarantees a graph to have a perfect matching.

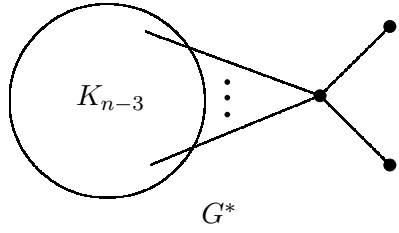
\begin{figure}[t]
	\setlength{\unitlength}{0.9pt}
	\begin{center}
		\begin{picture}(158.1,83.4)
		\qbezier(81.2,42.8)(81.2,26.0)(69.3,14.1)\qbezier(69.3,14.1)(57.4,2.2)(40.6,2.2)\qbezier(40.6,2.2)(23.8,2.2)(11.9,14.1)\qbezier(11.9,14.1)(0.0,26.0)(0.0,42.8)\qbezier(0.0,42.8)(0.0,59.6)(11.9,71.5)\qbezier(11.9,71.5)(23.8,83.4)(40.6,83.4)\qbezier(40.6,83.4)(57.4,83.4)(69.3,71.5)\qbezier(69.3,71.5)(81.2,59.6)(81.2,42.8)
		\put(129.1,45.0){\circle*{5}}
		\qbezier(56.6,71.8)(92.8,58.4)(129.1,45.0)
		\qbezier(129.1,45.0)(93.9,30.5)(58.7,16.0)
		\put(158.1,74.0){\circle*{5}}
		\qbezier(129.1,45.0)(143.6,59.5)(158.1,74.0)
		\put(158.1,16.0){\circle*{5}}
		\qbezier(129.1,45.0)(143.6,30.5)(158.1,16.0)
		\put(27.6,49.3){\makebox(0,0)[tl]{$K_{n-3}$}}
		\put(79.8,0.0){\makebox(0,0)[tl]{$G^*$}}
		\put(90.5,50.7){\circle*{2}}
		\put(90.5,37.1){\circle*{2}}
		\put(90.5,44.2){\circle*{2}}
		\end{picture}
	\end{center}
	\caption{The extremal graph $G^*$ of Theorem \ref{pm}.}
	\label{gstar}
\end{figure}

\begin{thm}\label{pm}
	Let $G$ be a connected graph with order $n$ and $n\ge 4$ be an even integer.
	\begin{enumerate}[(i)]
		\item For $ n \le 10$, if ${\lambda }_{1} (D\left(G\right))\le {\lambda }_{1} (D(S_{n,{\frac{n}{2}}-1}))$, then $G$ contains a perfect matching unless $G\cong S_{n,{\frac{n}{2}-1}}$.
		\item For $n\ge 12$, if ${\lambda }_{1} (D\left(G\right))\le {\lambda }_{1} (D(G^*))$, then $G$ contains a perfect matching unless $G\cong G^*$ where $G^*\cong K_1\vee (K_{n-3}\cup2K_1)$ (see Fig. \ref{gstar}).
	\end{enumerate}
\end{thm}

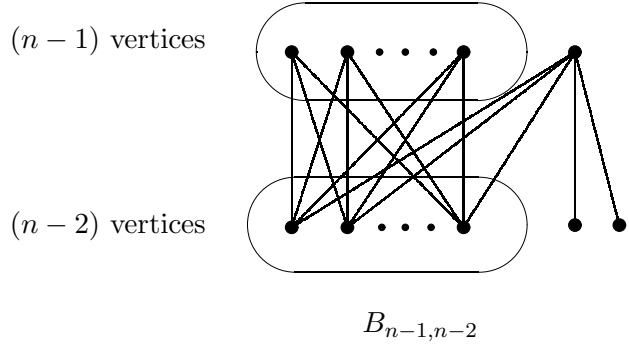
\begin{figure}[t]
	\setlength{\unitlength}{1.2pt}
	\begin{center}
		\begin{picture}(150.1,92.1)
		\put(77.6,22.5){\oval(87.0,29.7)}\put(136.3,76.9){\circle*{4}}
		\put(101.5,76.9){\circle*{4}}
		\put(65.3,76.9){\circle*{4}}
		\put(47.9,21.8){\circle*{4}}
		\put(101.5,21.8){\circle*{4}}
		\put(65.3,21.8){\circle*{4}}
		\put(76.1,21.8){\circle*{2}}
		\put(91.4,21.8){\circle*{2}}
		\put(83.4,21.8){\circle*{2}}
		\put(75.4,76.9){\circle*{2}}
		\put(92.1,76.9){\circle*{2}}
		\put(83.4,76.9){\circle*{2}}
		\qbezier(136.3,76.9)(92.1,49.3)(47.9,21.8)
		\qbezier(65.3,76.9)(65.3,49.3)(65.3,21.8)
		\qbezier(101.5,76.9)(101.5,49.3)(101.5,21.8)
		\qbezier(136.3,76.9)(100.8,49.3)(65.3,21.8)
		\qbezier(136.3,76.9)(118.9,49.3)(101.5,21.8)
		\qbezier(65.3,76.9)(56.6,49.3)(47.9,21.8)
		\qbezier(65.3,76.9)(83.4,49.3)(101.5,21.8)
		\qbezier(101.5,76.9)(74.7,49.3)(47.9,21.8)
		\qbezier(101.5,76.9)(83.4,49.3)(65.3,21.8)
		\put(150.1,22.5){\circle*{4}}
		\qbezier(136.3,76.9)(143.2,49.7)(150.1,22.5)
		\put(136.3,22.5){\circle*{4}}
		\qbezier(136.3,76.9)(136.3,49.7)(136.3,22.5)
		\put(-40,85.0){\makebox(0,0)[tl]{$(n-1)$ vertices}}
		\put(-40,27.0){\makebox(0,0)[tl]{$(n-2)$ vertices}}
		\put(70.0,-5.0){\makebox(0,0)[tl]{$B_{n-1,n-2}$}}
		\put(79.0,76.9){\oval(84.1,30.5)}\put(47.9,76.9){\circle*{4}}
		\qbezier(47.9,76.9)(47.9,49.3)(47.9,21.8)
		\qbezier(47.9,76.9)(56.6,49.3)(65.3,21.8)
		\qbezier(47.9,76.9)(74.7,49.3)(101.5,21.8)
		\end{picture}
	\end{center}
	
	\caption{The extremal graph $B_{n-1,n-2}$ of Theorem \ref{bppm}.}
	\label{bstar}
\end{figure}

A connected bipartite graph with two parts of size $ n_1 $ and $ n_2 $. We say that it
is balanced if $  n_1 = n_2 $. Note that if a connected bipartite graph is not balanced. Then it has no perfect matching. Let $B_{n-1 ,n-2}$ be the graph obtained from $K_{n,n-2}$ by attaching two pendent vertices to a vertex in $n$-vertex part (see Fig. \ref{bstar}). Then we have the following result.

\begin{thm}\label{bppm}
Let $G$ be a connected balanced bipartite graph with order $ 2n $ where $ n $ is an integer and $ n\ge 3 $. 
If $\lambda_{1}(D(G))\le \lambda_{1}(D(B_{n-1,n-2})) $, then $G$ has a perfect matching unless $G\cong B_{n-1,n-2}$ (see Fig. \ref{bstar}).
\end{thm}

\section{Proofs}

For $S\subseteq V(G) $, the subgraph induced by $V(G)-S$ is denoted by $G-S$.
A component is called an odd (even) component if the number of vertices in this component is odd (even) and let $o(G)$ denote the number of odd components of $ G $.
\begin{figure}[t]
	\setlength{\unitlength}{1pt}
	\begin{center}
		\begin{picture}(258.1,115.3)
		\qbezier(175.6,56.6)(175.6,46.9)(168.8,40.0)\qbezier(168.8,40.0)(161.9,33.2)(152.3,33.2)\qbezier(152.3,33.2)(142.6,33.2)(135.7,40.0)\qbezier(135.7,40.0)(128.9,46.9)(128.9,56.5)\qbezier(128.9,56.5)(128.9,66.2)(135.7,73.1)\qbezier(135.7,73.1)(142.6,79.9)(152.2,79.9)\qbezier(152.2,79.9)(161.9,79.9)(168.8,73.1)\qbezier(168.8,73.1)(175.6,66.2)(175.6,56.6)
		\qbezier(114.6,64.5)(114.6,58.8)(108.5,54.8)\qbezier(108.5,54.8)(102.4,50.8)(93.9,50.8)\qbezier(93.9,50.8)(85.3,50.8)(79.3,54.8)\qbezier(79.3,54.8)(73.2,58.8)(73.2,64.5)\qbezier(73.2,64.5)(73.2,70.2)(79.3,74.3)\qbezier(79.3,74.3)(85.3,78.3)(93.9,78.3)\qbezier(93.9,78.3)(102.4,78.3)(108.5,74.3)\qbezier(108.5,74.3)(114.6,70.2)(114.6,64.5)
		\qbezier(230.6,101.5)(230.6,95.8)(224.3,91.8)\qbezier(224.3,91.8)(218.0,87.7)(209.2,87.7)\qbezier(209.2,87.7)(200.3,87.7)(194.0,91.8)\qbezier(194.0,91.8)(187.8,95.8)(187.8,101.5)\qbezier(187.8,101.5)(187.8,107.2)(194.0,111.2)\qbezier(194.0,111.2)(200.3,115.3)(209.2,115.3)\qbezier(209.2,115.3)(218.0,115.3)(224.3,111.2)\qbezier(224.3,111.2)(230.6,107.2)(230.6,101.5)
		\qbezier(230.6,65.6)(230.6,59.8)(224.3,55.6)\qbezier(224.3,55.6)(218.0,51.5)(209.2,51.5)\qbezier(209.2,51.5)(200.3,51.5)(194.0,55.6)\qbezier(194.0,55.6)(187.8,59.8)(187.8,65.6)\qbezier(187.8,65.6)(187.8,71.5)(194.0,75.6)\qbezier(194.0,75.6)(200.3,79.8)(209.2,79.8)\qbezier(209.2,79.8)(218.0,79.8)(224.3,75.6)\qbezier(224.3,75.6)(230.6,71.5)(230.6,65.6)
		\qbezier(114.6,100.4)(114.6,94.6)(108.5,90.4)\qbezier(108.5,90.4)(102.4,86.3)(93.9,86.3)\qbezier(93.9,86.3)(85.3,86.3)(79.3,90.4)\qbezier(79.3,90.4)(73.2,94.6)(73.2,100.4)\qbezier(73.2,100.4)(73.2,106.3)(79.3,110.4)\qbezier(79.3,110.4)(85.3,114.6)(93.9,114.6)\qbezier(93.9,114.6)(102.4,114.6)(108.5,110.4)\qbezier(108.5,110.4)(114.6,106.3)(114.6,100.4)
		\qbezier(114.6,14.1)(114.6,8.3)(108.5,4.1)\qbezier(108.5,4.1)(102.4,0.0)(93.9,0.0)\qbezier(93.9,0.0)(85.3,0.0)(79.3,4.1)\qbezier(79.3,4.1)(73.2,8.3)(73.2,14.1)\qbezier(73.2,14.1)(73.2,20.0)(79.3,24.1)\qbezier(79.3,24.1)(85.3,28.3)(93.9,28.3)\qbezier(93.9,28.3)(102.4,28.3)(108.5,24.1)\qbezier(108.5,24.1)(114.6,20.0)(114.6,14.1)
		\qbezier(232.0,14.1)(232.0,8.3)(225.7,4.1)\qbezier(225.7,4.1)(219.5,0.0)(210.6,0.0)\qbezier(210.6,0.0)(201.8,0.0)(195.5,4.1)\qbezier(195.5,4.1)(189.2,8.3)(189.2,14.1)\qbezier(189.2,14.1)(189.2,20.0)(195.5,24.1)\qbezier(195.5,24.1)(201.8,28.3)(210.6,28.3)\qbezier(210.6,28.3)(219.5,28.3)(225.7,24.1)\qbezier(225.7,24.1)(232.0,20.0)(232.0,14.1)
		\qbezier(107.3,99.3)(122.9,83.7)(138.5,68.2)
		\qbezier(106.6,64.5)(120.7,61.6)(134.9,58.7)
		\qbezier(104.4,14.5)(121.1,31.2)(137.8,47.9)
		\qbezier(197.2,101.5)(180.9,85.2)(164.6,68.9)
		\qbezier(195.8,64.5)(182.3,61.3)(168.9,58.0)
		\qbezier(197.9,13.8)(181.3,30.5)(164.6,47.1)
		\put(95.0,45.7){\circle*{2}}
		\put(95.0,32.6){\circle*{2}}
		\put(95.0,39.2){\circle*{2}}
		\put(210.3,46.4){\circle*{2}}
		\put(210.3,34.1){\circle*{2}}
		\put(210.3,40.6){\circle*{2}}
		\put(87.7,105.1){\makebox(0,0)[tl]{$G_1$}}
		\put(87.7,68.9){\makebox(0,0)[tl]{$G_2$}}
		\put(87.7,18.9){\makebox(0,0)[tl]{$G_q$}}
		\put(202.3,106.6){\makebox(0,0)[tl]{$R_1$}}
		\put(202.3,69.6){\makebox(0,0)[tl]{$R_2$}}
		\put(205.2,18.9){\makebox(0,0)[tl]{$R_k$}}
		\put(149.2,61.1){\makebox(0,0)[tl]{$S$}}
		\put(147.2,-7.7){\makebox(0,0)[tl]{$G$}}
		\put(27.6,63.8){\makebox(0,0)[tl]{odd}}
		\put(0.0,49.3){\makebox(0,0)[tl]{components}}
		\put(245.1,66.7){\makebox(0,0)[tl]{even}}
		\put(229.8,54.4){\makebox(0,0)[tl]{components}}
		\end{picture}
	\end{center}
	\caption{The graph in Tutte's theorem.}
	\label{ttgraph}
\end{figure}
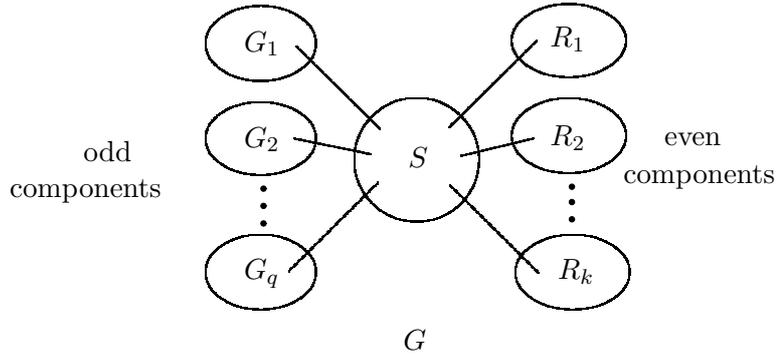

\begin{lemma}[Tutte's theorem \cite{Tutte}]\label{tt}
A graph $ G $ has a perfect matching (see Fig.\ref{ttgraph}) if and only if\par
\centering $ o(G- S) \le  |S| $ for each $ S \subseteq V(G) $.
\end{lemma}

Let $W(G)=\sum_{i<j}d_{ij}$ be the \emph{Wiener index} of a connected graph $G$ with order $n$.
Note that $\lambda_1(D(G))=\max \limits_{\textbf{x}\in {\mathbb{R}}^{n}}\frac{\textbf{x}^{t}D(G)\textbf{x}}{\textbf{x}^{t}\textbf{x}}$.
Then we have
$$\lambda_1(D(G))=\max \limits_{\textbf{x}\in {\mathbb{R}}^{n}}\frac{\textbf{x}^{t}D(G)\textbf{x}}{\textbf{x}^{t}\textbf{x}}\geq \frac{\mathbf{1}^{t}D\mathbf{1}}{\mathbf{1}^{t}\mathbf{1}}\geq \frac{2W(G)}{n},$$
where $\mathbf{1}=(1,1,\ldots,1)^{t}$. Now we give the proof of Theorem \ref{pm}.
{\flushleft \textbf{Proof of Theorem \ref{pm}.}}
By way of contradiction assume that $ G $ has no perfect matching with the minimum distance spectral radius. By Lemma \ref{tt}, there exists $ S \subset V (G) $ such that $  q -|S| > 0 $, where $ o(G-S) = q $ and all components of $ G-S $ are odd, otherwise, we can remove one vertex from each even component to the set $S$, in consequence, the number of odd component and the number of vertices in $S$ have the same increase, so that $q$ is always larger than $|S|$ and all components of $G-S $ are odd. Note that $G$ is connected, which implies that $S$ is not empty. $|S|$ and $q$ have the same parity because $ n $ is even, so $q-|S| \ge 2$. Let $ G'$ be the graph obtained from $ G $ by joining $ S $ and $ G-S $ and by adding edges in $ S $ and in all components of $ G-S $ so that all components of $ G-S $ and $ G[S] $ are cliques. It's clear that ${G}^{\prime }\cong {K}_{s}\vee \left({K}_{n_1}\cup K_{n_2}\cup \cdots \cup {K}_{n_q}\right)$ where $ |S|=s$, ${n}_{1}\ge {n}_{2}\ge \dots \ge {n}_{q}\ge 1$ and ${n}_{1}+{n}_{2}+ \dots + {n}_{q}= n-s$. Note that $G'\cong G$, otherwise,
referring to the Perron-Frobenius theorem, we have ${\lambda }_{1}\left(D\left(G'\right)\right)<{\lambda }_{1}\left(D\left(G\right)\right)$, a contradiction. Suppose ${n}_{2}={n}_{3}= \dots = {n}_{q}= 1$, then we get a graph ${G}^{\prime \prime }\cong {K}_{s}\vee \left({K}_{n-s-\left(q-1\right)}\cup (q-1)K_1\right)$.

\begin{claim}
	${\lambda }_{1}\left(D\left(G'\right)\right)\ge {\lambda }_{1}\left(D\left(G''\right)\right)$ with equality if and only if $G' \cong G''$.
\end{claim}
If $n_1=1$, then $n_1={n}_{2}={n}_{3}= \dots = {n}_{q}= 1$ and $G' \cong G''$. Now we consider $n_1\ge 3$.
Denote the vertex set of $G''$ by $V(G'')=V(K_{s})\cup V(K_{n-s-(q-1)}) \cup V((q-1)K_{1})$.
Suppose that $X$ is the Perron vector of $D(G'')$,
and let $x(v)$ denote the entry of $X$ corresponding to the vertex $v\in V(G'')$.
By symmetry, it is easy to see that all vertices of $V(K_{s})$ (resp. $V(K_{n-s-(q-1)})$ and $V((q-1)K_{1})$) have the same entries in $X$. Thus we can suppose $x(u)=a$ for any $u\in V((q-1)K_{1})$, $x(v)=b$ for any $v\in V(K_{n-s-(q-1)})$ and $x(w)=c$ for any  $w\in V(K_{s})$.
Then
\[
\left\{
\begin{array}{l}
\lambda_1(D(G''))a=sc+2\left(n-s-q+1\right)b+2\left(q-2\right)a,\\[3mm]
\lambda_1(D(G''))b=sc+\left(n-s-q\right)b+2\left(q-1\right)a,\\[3mm]
\lambda_1(D(G''))c=(s-1)c+\left(n-s-q+1\right)b+\left(q-1\right)a.
\end{array}
\right.
\]
Thus,
$$a=\left[1+\frac{n-s-q}{\lambda_1(D(G''))+2}\right]b.$$ It follows that
\begin{eqnarray*}
	&&{\lambda }_{1}\left(D\left(G'\right)\right)- {\lambda }_{1}\left(D\left(G''\right)\right)~~~~~~~~~~~~~~~~~~~~~~~~~~~~~~~~~~~~~~~~~~~~~~~~~~~~~~~~~~~~~~~~~~~~~~~~~~~~~~~~~\\[2mm]
	&\ge& {X}^{t}\left(D\left({G}^{\prime }\right)-D\left({G}^{\prime \prime }\right)\right)X\\
	&=&n_1\sum _{k=2}^{q}\left({n}_{k}-1\right){b}^{2}
	+\left({n}_{2}-1\right)\left[\left(n-s-{n}_{2}-\left(q-2\right)\right){b}^{2}-2ab\right]\\
	&+&\left({n}_{3}-1\right)\left[\left(n-s-{n}_{3}-\left(q-2\right)\right){b}^{2}-2ab\right]+\left({n}_{4}-1\right)\left[\left(n-s-{n}_{4}-\left(q-2\right)\right){b}^{2}-2ab\right]\\
	&+&\cdots +\left({n}_{q}-1\right)\left[\left(n-s-{n}_{q}-\left(q-2\right)\right){b}^{2}-2ab\right].
\end{eqnarray*}
Since $ n_1\ge 3 $ and $ n_2\ge n_3 \ge \cdots\ge n_q\ge 1 $, in order to show ${\lambda }_{1}\left(D\left(G'\right)\right)- {\lambda }_{1}\left(D\left(G''\right)\right)>0$, we only need to prove $\left(n-s-{n}_{2}-\left(q-2\right)\right){b}^{2}-2ab>0$.

Note that $K_{n-q+1}$ is a subgraph of $G''$. Then ${\lambda }_{1}\left(D\left(G''\right)\right) >{\lambda }_{1}\left(D\left(K_{n-q+1}\right)\right)=n-q,$ which implies
\begin{eqnarray*}
	&&\left(n-s-{n}_{2}-\left(q-2\right)\right){b}^{2}-2ab~~~~~~~~~~~~~~~~~~~~~~~~~~~~~~~~~~~~~~~~~~~~~~~~~~~~~~~~~~~~~~~~~~~~~~~~~~~~~~~\\
	&=&b^2\left(n-s-n_2-q-\frac{2n-2s-2q}{\lambda_1(D(G''))+2}\right)\\[3mm]
	&>&b^2\left(n-s-n_2-q-\frac{2n-2s-2q}{n-q+2}\right)\\[3mm]
	&=&b^2\left(n-s-n_2-q-2+\frac{2s+4}{n-q+2}\right)\\[3mm]
	&>&{b}^{2}\left(n-s-{n}_{2}-q-2\right)\\[3mm]
	&=&{b}^{2}(n_1+n_2+\cdots+n_{q-1}+n_q-n_2-q-2)\\[3mm]
	&=&{b}^{2}\left(\sum _{i=1, i\ne 2}^{q}\left({n}_{i}-1\right)-1\right)\\&>& 0.
\end{eqnarray*}
So Claim 1 holds.

So $G\cong G''$, or there will be a contradiction. We know that $q-|S|\ge 2$. Let $\widetilde{G}\cong {K}_{s}\vee \left({K}_{n-2s-1}\cup \left(s+1\right){K}_{1}\right)$. We compare the distance spectral radius of $G''$ and $\widetilde{G}$ in the following.
\begin{claim}
	${\lambda }_{1}\left(D\left(G''\right)\right)\ge {\lambda }_{1}(D(\widetilde{G}))$ with equality if and only if $G'' \cong \widetilde{G}$.
\end{claim}
Recall ${G}^{\prime \prime }\cong {K}_{s}\vee \left({K}_{n-s-\left(q-1\right)}\cup (q-1)K_1\right)$. If $q=s+2$, then $\widetilde{G} \cong G''$. Now we consider $q\ge s+4$.
Denote the vertex set of $\widetilde{G}$ by $V(\widetilde{G})=V(K_{s})\cup V(K_{n-2s-1}) \cup V((s+1)K_{1})$.
Suppose that $Y$ is the Perron vector of $D(\widetilde{G})$,
and let $Y(v)$ denote the entry of $Y$ corresponding to the vertex $v\in V(\widetilde{G})$.
By symmetry, it is easy to see that all vertices of $V(K_{s})$ (resp. $V(K_{n-2s-1})$ and $V((s+1)K_{1})$) have the same entries in $X$. Thus we can suppose $Y(u)=y_1$ for any $u\in V((s+1)K_{1})$, $Y(v)=y_2$ for any $v\in V(K_{n-2s-1})$ and $Y(w)=y_3$ for any  $w\in V(K_{s})$. Let $n_1'=n-s-(q-1)\ge 1$ and ${G}^{\prime \prime }\cong {K}_{s}\vee \left({K}_{n_1'}\cup (q-1)K_1\right)$. Then
\begin{eqnarray*}
	&&{\lambda }_{1}\left(D\left(G''\right)\right)- {\lambda }_{1}(D(\widetilde{G}))~~~~~~~~~~~~~~~~~~~~~~~~~~~~~~~~~~~~~~~~~~~~~~~~~~~~~~~~~~~~~~~~~~~~~~~~~~~~~~~~~\\
	&\ge& {Y}^{t}\left(D\left({G''}\right)-D(\widetilde{G})\right)Y\\
	&=&n_1'(q-s-2)y_2^2+(n_1'+q-s-3)(q-s-2)y_2^2\\
	&=&y_2^2\left[{q}^{2}+\left(2{n}_{1}^{\prime }-2s-5\right)q+{s}^{2}+5s-2{n}_{1}^{\prime }s-4{n}_{1}^{\prime }+6\right].
\end{eqnarray*}
Since $$\frac{-\left(2{n}_{1}^{\prime }-2s-5\right)}{2}=-{n}_{1}^{\prime }+s+\frac{5}{2}<s+4,$$
we obtain
\begin{eqnarray*}
	&&y_2^2\left[{q}^{2}+\left(2{n}_{1}^{\prime }-2s-5\right)q+{s}^{2}+5s-2{n}_{1}^{\prime }s-4{n}_{1}^{\prime }+6\right]~~~~~~~~~~~~~~~~~~~~~~~~~~~~~~~~~~~~~~~~~~~~~~~~~~~~~~~~~~~~~~~~~\\
	&\ge &y_2^2\left[{(s+4)}^{2}+\left(2{n}_{1}^{\prime }-2s-5\right)(s+4)+{s}^{2}+5s-2{n}_{1}^{\prime }s-4{n}_{1}^{\prime }+6\right]\\
	&=&y_2^2(4{n}_{1}^{\prime }+2)\\&>& 0.
\end{eqnarray*}
So Claim 2 holds.

So $G\cong \widetilde{G}$, or there will be a contradiction. Let $  G^* \cong K_1\vee (K_{n-3}\cup 2K_1)$. In the end, we shall show that $G^*$ contains the minimum distance spectral radius and no perfect matching under the most situations.

\begin{claim}
	If $ n\ge 2s+4 $, then ${\lambda }_{1}(D(\widetilde{G}))\ge {\lambda }_{1}\left(D\left(G^*\right)\right)$ with equality if and only if $\widetilde{G} \cong G^*$.
\end{claim}
Recall $\widetilde{G}\cong {K}_{s}\vee \left({K}_{n-2s-1}\cup \left(s+1\right){K}_{1}\right)$. If  $s=1$, then $G^*\cong \widetilde{G}$. Now we suppose $s\ge 2$ so that $n\ge 2s+4\ge 8$. Then the quotient matrix of the partition $\{V(K_{n-2s-1},V(K_s),V((s+1)K_1))\}$ of $\widetilde{G}$ is

$$\left(\begin{array}{ccccccc}
n-2s-2 & s & 2(s+1)\\
\ &\ &\ \\
n-2s-1 & s-1 & s+1\\
\ &\ &\ \\
2(n-2s-1) & s & 2s
\end{array}\right),$$
the characteristic polynomial of the matrix is
\begin{eqnarray*}
	f(x)={x}^{3}-\left(s+n-3\right){x}^{2}-\left(2sn-5{s}^{2}+5n-6s-6\right)x+n{s}^{2}-2{s}^{3}-sn+2{s}^{2}-4n+6s+4.&&~~~~~~~~~~~~~~~~~~~~~~~~~~~~
\end{eqnarray*}
We know that $ {\lambda }_{1}(D(\widetilde{G}))$ is the largest root of $ f(x)=0 $. Since $ {\lambda }_{1}\left(D\left(G^*\right)\right)= \theta(n) $ (simply $\theta$) is the largest root of the equation $q(x)={x}^{3}-\left(n-2\right){x}^{2}-\left(7n-17\right)x-4n+10=0$, we have
\begin{eqnarray*}
	h(\theta)&=&f\left(\theta \right)-q\left(\theta \right)\\
	&=&-\left(s-1\right){\theta }^{2}-\left(2ns-5{s}^{2}-2n-6s+11\right)\theta +n{s}^{2}-2{s}^{3}-sn+2{s}^{2}+6s-6\\
	&=&(s - 1)(-\theta^2 + (-2n + 5s + 11)\theta + sn - 2s^2 + 6).
\end{eqnarray*}
Moreover, $\theta={\lambda }_{1}\left(D\left(G^*\right)\right)\ge \frac{2W(G^*)}{n}=\frac{n^2+3n-10}{n}\ge n+1$ and $s\ge 2$.
To show ${\lambda }_{1}(D(\widetilde{G}))>{\lambda }_{1}\left(D\left(G^*\right)\right)$, we only need to prove $ h_1(\theta) =-\theta^2 + (-2n + 5s + 11)\theta + sn - 2s^2 + 6< 0$ when $\theta\ge n+1$, then $h(\theta)<h_1(\theta)<0$. Note that $$\frac{-2n+5s+11}{2}=-n+\frac{5}{2}s+\frac{11}{2}< n+1.$$  So when $\theta\ge n+1$, $ h_1(\theta)$  monotonically decreases as $\theta$ increases and $ h_1(\theta)\le h_1(n+1) $.
Let $ g(n)=h_1(n+1)=-3n^2 + (6s + 7)n - 2s^2 + 5s + 16 $ where $ n\ge 2s+4 $.
Then
$$\frac{6{s}+7}{6}=s+\frac{7}{6}<2s+4.$$ So when $n\ge 2s+4$, $ g(n)$ is monotonically decreasing and $ h(\theta)<h_1(\theta)\le g(n)\le  g(2s+4)=-2s^2 - 5s - 4<0$. So Claim 3 holds.

%Therefore, $G\cong G^*$ if $G$ has no perfect matching with the minimum distance spectral radius, which is contrary to $G\ncong G^*$.

If $ n=2s+2 $, observe that $\widetilde{G}\cong S_{n,{\frac{n}{2}-1}}$ and the quotient matrix of the partition $\{V(K_{\frac{n}{2}-1}),V(({\frac{n}{2}+1})K_1)\}$ of $S_{n,{\frac{n}{2}-1}}$ is
$$\left(\begin{array}{cccccccc}
{\frac{n}{2}-2} & \ &{\frac{n}{2}+1}\\
\ &\ &\ \\
{\frac{n}{2}-1}&\ & n
\end{array}\right).$$
By a simple calculation, $$\lambda_{1}(D( S_{n,{\frac{n}{2}-1}}))=\frac{3n-4+\sqrt{n^2+24n-16}}{4}.$$

\begin{claim}
	If $ n=2s+2 $ and $n\ge 12$, then ${\lambda }_{1}(D(S_{n,{\frac{n}{2}-1}}))> {\lambda }_{1}\left(D\left(G^*\right)\right)$.
\end{claim}
Let $n=2s+2$, then  $g(n)<0$ when $ n\ge 12 $. So Claim 4 holds.

\begin{claim}
	If $ n=2s+2 $ and $n\le 10$, then ${\lambda }_{1}(D(S_{n,{\frac{n}{2}-1}}))\le {\lambda }_{1}\left(D\left(G^*\right)\right)$ with equality if and only if $G^*\cong S_{n,{\frac{n}{2}-1}}$.
\end{claim}

If $ s=1$ and $n=4 $, then $S_{n,{\frac{n}{2}-1}}\cong G^*$.
If $ s=2 $ and $ n=6 $, then $\lambda_{1}(D( S_{6,2}))=\frac{7+\sqrt{41}}{2}$ and $ q(\frac{7+\sqrt{41}}{2})<0 $, so ${\lambda }_{1}(D(S_{6,2}))< \theta(6)$.
If $ s=3 $ and $ n=8 $, then $\lambda_{1}(D( S_{8,3}))=5+\sqrt{15}$ and $ q(5+\sqrt{15})<0 $, so 	${\lambda }_{1}(D(S_{8,3}))< \theta(8)$.
If $ s=4$ and $ n=10 $, then $\lambda_{1}(D( S_{10,4}))=11$ and $ q(11)<0 $, so 	${\lambda }_{1}(D(S_{10,4}))< \theta(10)$. So Claim 5 holds.

In conclusion, $G\cong G^*$ if $n\ge 12$ and $G\cong S_{n,{\frac{n}{2}-1}}$ if $4\le n\le 10 $, otherwise $G$ doesn't contain the minimum distance spectral radius, a contradiction.
\qed

In some specific applications, one needs to find a matching in a bipartite graph which can cover one partite. 
Essential and sufficient conditions for the existence of such a matching was first proposed by Hall \cite{hallthm}. Let $ S $ be a set
of vertices in a graph $ G $. The set of all neighbours of the vertices in $ S $ is denoted
by $ N(S) $.

\begin{lemma}[Hall's theorem \cite{hallthm}] \label{hall}
	A bipartite graph $ G := G[X,Y ] $ has a matching which covers every vertex in $ X $ if
	and only if\par
	\centering $  |N(S)| \ge |S|  $ for each $  S \subseteq X.  $
\end{lemma}

%\begin{thm}\label{thm4}
%Suppose that $G$ is a connected nearly-balanced bipartite graph of order $2n+1$. If $\lambda_1(D(G))<\lambda_1(D(B_n^2))$,
%then $G$ contains a $P_{2n+1}$.
%\end{thm}

\begin{figure}[t]
	\setlength{\unitlength}{1pt}
	\begin{center}
		\begin{picture}(190.0,140.7)
		\put(43.5,110.2){\oval(87.7,29.7)}\put(146.5,109.5){\oval(86.3,29.0)}\put(43.5,48.6){\oval(87.7,31.2)}\put(146.5,49.3){\oval(87.0,29.7)}\put(125.4,47.9){\circle*{4}}
		\put(125.4,109.5){\circle*{4}}
		\put(13.1,108.8){\circle*{4}}
		\put(43.5,108.8){\circle*{4}}
		\put(28.3,108.8){\circle*{4}}
		\put(166.8,48.6){\circle*{4}}
		\put(166.8,109.5){\circle*{4}}
		\put(18.9,47.9){\circle*{4}}
		\put(43.5,47.9){\circle*{4}}
		\put(76.1,108.8){\circle*{4}}
		\put(76.1,47.9){\circle*{4}}
		\qbezier(13.1,108.8)(16.0,78.3)(18.9,47.9)
		\qbezier(28.3,108.8)(23.6,78.3)(18.9,47.9)
		\qbezier(43.5,108.8)(31.2,78.3)(18.9,47.9)
		\qbezier(18.9,47.9)(47.5,78.3)(76.1,108.8)
		\qbezier(13.1,108.8)(28.3,78.3)(43.5,47.9)
		\qbezier(28.3,108.8)(35.9,78.3)(43.5,47.9)
		\qbezier(43.5,108.8)(43.5,78.3)(43.5,47.9)
		\qbezier(76.1,108.8)(59.8,78.3)(43.5,47.9)
		\qbezier(13.1,108.8)(44.6,78.3)(76.1,47.9)
		\qbezier(28.3,108.8)(52.2,78.3)(76.1,47.9)
		\qbezier(43.5,108.8)(59.8,78.3)(76.1,47.9)
		\qbezier(76.1,108.8)(76.1,78.3)(76.1,47.9)
		\qbezier(18.9,47.9)(72.1,78.7)(125.4,109.5)
		\qbezier(18.9,47.9)(92.8,78.7)(166.8,109.5)
		\qbezier(43.5,47.9)(84.5,78.7)(125.4,109.5)
		\qbezier(43.5,47.9)(105.1,78.7)(166.8,109.5)
		\qbezier(125.4,109.5)(100.8,78.7)(76.1,47.9)
		\qbezier(76.1,47.9)(121.4,78.7)(166.8,109.5)
		\qbezier(125.4,109.5)(125.4,78.7)(125.4,47.9)
		\qbezier(166.8,109.5)(146.1,78.7)(125.4,47.9)
		\qbezier(125.4,109.5)(146.1,79.0)(166.8,48.6)
		\qbezier(166.8,109.5)(166.8,79.0)(166.8,48.6)
		\put(136.3,109.5){\circle*{2}}
		\put(155.2,109.5){\circle*{2}}
		\put(145.7,109.5){\circle*{2}}
		\put(52.2,47.9){\circle*{2}}
		\put(68.2,47.9){\circle*{2}}
		\put(60.2,47.9){\circle*{2}}
		\put(137.0,47.9){\circle*{2}}
		\put(155.2,47.9){\circle*{2}}
		\put(145.7,47.9){\circle*{2}}
		\put(52.9,108.8){\circle*{2}}
		\put(67.4,108.8){\circle*{2}}
		\put(60.2,108.8){\circle*{2}}
		\put(39.2,140.7){\makebox(0,0)[tl]{$S$}}
		\put(34.1,30.5){\makebox(0,0)[tl]{$N(S)$}}
		\put(129.8,139.9){\makebox(0,0)[tl]{$X-S$}}
		\put(120.4,31.9){\makebox(0,0)[tl]{$Y-N(S)$}}
		\put(74.7,0.0){\makebox(0,0)[tl]{$G[X,Y]$}}
		\end{picture}
	\end{center}
	\caption{The graph $B_{|S|,|N(S)|}$ in Theorem \ref{bppm}.}
	\label{hallgraph}
\end{figure}
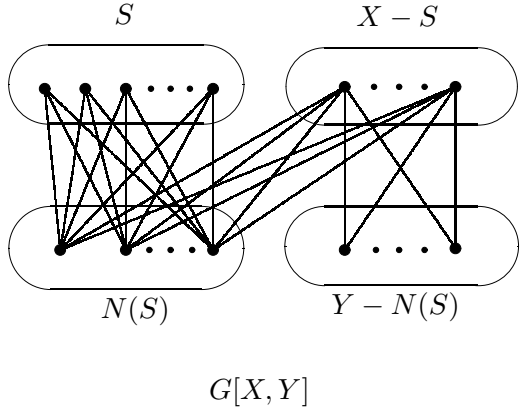

For two vertex sets $X$ and $Y$, let $e(X,Y)$ be the set of all edges between $X$ and $Y$.
Now we are in a position to prove Theorem \ref{bppm}.

{\flushleft \textbf{Proof of Theorem \ref{bppm}.}}
Assume to the contrary that $ G $ has no perfect matching with the minimum distance spectral radius. Let $ G:=G[X,Y]  $ be a connected balanced bipartite graph, where $ |X|=|Y|=n $. By Lemma \ref{bppm}, since $G$ has no perfect matching, there exists $ S\subset X $ and $ |N(S)|<|S| $. Notice that there exists no edges between $ S $ and $ Y-N(S)  $, otherwise, we can find a vertex $ v\in S $ and $ Y-N(S) $ contains its neighbors, a contradiction. Let $ B_{s,p} $ be the connected balanced bipartite graph obtained from $ G $ by joining $ S $ and $ N(S) $, $ X-S $ and $ Y-N(S) $  and  by adding all possible edges between $ X-S $ and $  N(S) $ where $|S|=s$, $|N(S)|=p $ and $ 1\le p<s\le n-1$, so that $ B_{s,p}\cong K_{s,p} \cup K_{n-s,n-p} +e(N(S),X-S)$, i.e. $ B_{s,p}\cong K_{n,n} -e(S,Y-N(S))$ (see Fig.\ref{hallgraph}). Note that $G\cong B_{s,p} $, otherwise, we have ${\lambda }_{1}\left(D\left(B_{s,p}\right)\right)<{\lambda }_{1}\left(D\left(G\right)\right)$ by Perron-Frobenius theorem, a contradiction.
\setcounter{claim}{0}
\begin{claim}
	$\lambda_{1}(D(B_{s,p}))\ge \lambda_{1}(D(B_{s,s-1}))$ with equality if and only if $B_{s,p} \cong B_{s,s-1}$.
\end{claim}
If $p=s-1$, then $B_{s,p} \cong B_{s,s-1}$. Now we consider $1\le p\le s-2$.
Denote the vertex set of $B_{s,s-1}$ by $V(B_{s,s-1})=S\cup (X-S) \cup N(S) \cup (Y-N(S))$ where $|S|=s$ and $|N(S)|=s-1$. Referring to the Perron-Frobenius theorem, suppose that a positive vector $Z$ is the Perron vector of $D(B_{s,s-1})$,
and let $Z(v)$ denote the entry of $Z$ corresponding to the vertex $v\in V(B_{s,s-1})$.
By symmetry, it is easy to see that all vertices of $S$ (resp. $X-S$, $N(S)$ and $Y-N(S)$) have the same entries in $Z$. Thus we can suppose $Z(u)=z_1$ for any $u\in S$, $Z(v)=z_2$ for any $v\in X-S$, $Z(w)=z_3$ for any  $w\in N(S)$ and $Z(z)=z_4$ for any  $z\in Y-N(S)$. Therefore, we have
\begin{eqnarray*}
	&&{\lambda }_{1}\left(D\left(B_{s,p}\right)\right)- {\lambda }_{1}(D(B_{s,s-1}))~~~~~~~~~~~~~~~~~~~~~~~~~~~~~~~~~~~~~~~~~~~~~~~~~~~~~~~~~~~~~~~~~~~~~~~~~~~~~~\\&\ge& {Z}^{t}\left(D\left({B_{s,p}}\right)-D(B_{s,s-1})\right)Z\\
	&=&4s(s-1-p)z_1z_3\\&>&0.
\end{eqnarray*}
So Claim 1 holds.

Therefore, $G\cong B_{s,s-1}$, or there will be a contradiction.
\begin{claim}
	$\lambda_{1}(D(B_{s,s-1}))\ge \lambda_{1}(D(B_{n-1,n-2}))$ with equality if and only if $B_{s,s-1} \cong B_{n-1,n-2}$.
\end{claim}
If $s=n-1$, then $B_{s,s-1} \cong B_{n-1,n-2}$. Now consider  $2\le s\le n-2$. The quotient matrix of the partition $\{S, X-S, N(S), Y-N(S)\}$ of $B_{s,s-1}$ where $|S|=s\le n-2$ and $|N(S)|=s-1$ is

$$\left(\begin{array}{ccccccccc}
2s-2 \ \ & s-1 \ \ & 2n-2s \ \ & 3n-3s+3 \\[3mm]

s \ \ & 2s-4 \ \ & n-s \ \ & 2n-2s+2 \\[3mm]

2s \ \ & s-1 \ \ & 2n-2s-2 \ \ & n-s+1 \\[3mm]

3s \ \ & 2s-2 \ \ & n-s \ \ & 2n-2s \\[3mm]
\end{array}\right),$$
the characteristic polynomial of the matrix is
\begin{eqnarray*}
f(\lambda)=\lambda^4 &+& (-4n + 8)\lambda^3 + (3n^2 - 8ns + 8s^2 - 24n - 8s + 24)\lambda^2 \\&+& (8n^2s - 8ns^2 + 12n^2 - 32ns + 40s^2 - 48n - 40s + 32)\lambda\\ &-& 12n^2s^2 + 24s^3n - 12s^4 + 28n^2s - 52ns^2 + 24s^3 + 12n^2 \\&-& 20sn + 36s^2 - 32n - 48s + 16.
\end{eqnarray*}
We know that $ {\lambda }_{1}(D(B_{s,s-1}))$ is the largest root of $ f(\lambda )=0 $. Since $ {\lambda }_{1}\left(D\left(B_{n-1,n-2}\right)\right)$, written by $\rho $, is the largest root of the equation
$$q(\lambda)=\lambda^4 + (-4n + 8)\lambda^3 + (3n^2 - 40n + 40)\lambda^2 + (28n^2 - 144n + 112)\lambda + 20n^2 - 88n + 64=0,$$
we have
\begin{eqnarray*}
	h(\rho)&=&f\left(\rho \right)-q\left(\rho \right)\\
	&=&( 8s^2-(8n+8)s + 16n  - 16)\rho^2 \\&+& (8sn^2 - 8ns^2 - 16n^2 - 32ns + 40s^2 + 96n - 40s - 80)\rho \\&-& 12n^2s^2 + 24s^3n - 12s^4 + 28sn^2 - 52ns^2 + 24s^3 - 8n^2 - 20sn + 36s^2 + 56n - 48s - 48.
\end{eqnarray*}
Now we need to state $h(\rho)<0$. Firstly, we give a lower bound on $\rho$,
$$\frac{2W(B_{n-1,n-2})}{2n}=2n+5-\frac{6}{n}<\rho.$$
By a simple computation, we obtain
$$ 8s^2-(8n+8)s + 16n  - 16\le 0 \ \mbox{for}\  2\le s\le n-2,$$
and
\begin{eqnarray*}
\frac{8sn^2 - 8ns^2 - 16n^2 - 32ns + 40s^2 + 96n - 40s - 80}{-2( 8s^2-(8n+8)s + 16n  - 16)}=\frac{n-5}{2}<2n<\rho.
\end{eqnarray*}
Thus, $h(\rho)$ is monotonically decreasing for $\rho> 2n$ and $h(\rho)<h(2n)$. Let $h(2n)$ be written by $ g(s) $ and
we only need to prove $ g(s)<0 $ for $2\le s\le n-2$ in the following steps. Based on Matlab programming, we have
\begin{eqnarray*}
	g(s)=-4(s - 2)(-s + n - 1)[- 3s^2 + (3s+3)n+19n  +4n^2+ 6].
\end{eqnarray*}
It is easy to see $s-2\ge 0$ and $-s + n - 1>0$. Define $  r(s)=- 3s^2 + (3s+3)n+19n  +4n^2+ 6 $, we can easily get for $2\le s\le n-2$ and $n\ge 3$,
\begin{eqnarray*}
&&r(s)=- 3s^2 + (3s+3)n+19n  +4n^2+ 6\\&>&\min \{r(2),\ r(n-2)\}~~~~~~~~~~~~~~~~~~~~~~~~~~~~~~~~~~~~~~~~~~~~~~~~~~~~~~~~~~~~~~~~~~~~~~~~~~~~~~~~~~~~~~~~\\
&=&\min \{4n^2+25n,\ 4n^2+28n-12\}\\
&>&0.
\end{eqnarray*}
Therefore, $h(\rho)<h(2n)=g(s)<0$.
So Claim 2 holds.

In conclusion, $G\cong B_{n-1,n-2}$ has the minimum distance spectral radius among all $2n$-vertex balanced bipartite graphs without a perfect matching.\qed

%\noindent\textbf{Acknowledgement}

\end{document}